\documentclass[11pt]{amsart}

\usepackage{hyperref}
\usepackage{amsmath,amssymb,amsfonts}
\usepackage{amsthm}
\usepackage[curve]{xypic}
\usepackage{enumerate}
\usepackage{charter}

\theoremstyle{plain}
\newtheorem{thm}{Theorem}[subsection]
\newtheorem{prop}[thm]{Proposition}
\newtheorem{cor}[thm]{Corollary}
\newtheorem{lem}[thm]{Lemma}

\theoremstyle{definition}
\newtheorem{dfn}[thm]{Definition}

\newtheorem*{notation}{Notation}

\newtheorem{eg}[thm]{Example}
\newtheorem{rmk}[thm]{Remark}

\newenvironment{varthm}[1]{%
\par\addvspace{8pt plus3pt minus2pt}%
\noindent{\bf {}#1.\hskip.5em}\itshape\ignorespaces}{%
	\par\vskip-\baselineskip
        \par
	\addvspace{8pt plus3pt minus3pt}}

\newcommand{\Nguyen}{Nguy$\tilde{\hat{\text{e}}}$n}

\newcommand{\PP}{\mathbb{P}} 
\newcommand{\CC}{\mathbb{C}} 
\newcommand{\OO}{\mathcal{O}} 
\newcommand{\FF}{\mathbb{F}} 
\DeclareMathOperator{\Pic}{Pic} 
\DeclareMathOperator{\Hom}{Hom} 
\newcommand{\SU}{\mathcal{SU}} 
\newcommand{\U}{\mathcal{U}} 
\DeclareMathOperator{\gr}{gr} 
\newcommand{\calC}{\mathcal{C}} 
\newcommand{\calD}{\mathcal{D}} 
\newcommand{\calI}{\mathcal{I}} 
\newcommand{\calS}{\mathcal{S}} 
\DeclareMathOperator{\Fix}{Fix} 
\DeclareMathOperator{\Sym}{Sym} 
\DeclareMathOperator{\Sing}{Sing} 
\newcommand{\st}{\, : \,} 
\newcommand{\Thgen}{\Theta^{\text{gen}}} 
\newcommand{\bino}[2]{\bigl( \begin{smallmatrix} #1 \\ #2 \end{smallmatrix} \bigr)}

\begin{document}


\bibliographystyle{amsalpha}

\title{Vector bundles, dualities, and classical geometry on a curve of genus two}

\author[\Nguyen{} Q. M.]{\Nguyen{} Quang Minh}

\address{3065 Jackson St \\ San Francisco, CA 94115 (USA)}

\email{{\tt quang-minh.nguyen@sfuhs.org}}

\date{} 
 
\begin{abstract}
  Let $C$ be a curve of genus two.  We denote by $\SU_C(3)$ the moduli space of semi-stable vector
  bundles of rank 3 and trivial determinant over $C$, and by $J^d$ the variety of line bundles of
  degree $d$ on $C$.  In particular, $J^1$ has a canonical theta divisor $\Theta$.  The space
  $\SU_C(3)$ is a double cover of $\PP^8=|3\Theta|$ branched along a sextic hypersurface, the Coble
  sextic.  In the dual $\check{\PP}^8=|3\Theta|^*$, where $J^1$ is embedded, there is a unique cubic
  hypersurface singular along $J^1$, the Coble cubic.  We prove that these two hypersurfaces are
  dual, inducing a non-abelian Torelli result.  Moreover, by looking at some special linear sections
  of these hypersurfaces, we can observe and reinterpret some classical results of algebraic
  geometry in a context of vector bundles: the duality of the Segre-Igusa quartic with the Segre
  cubic, the symmetric configuration of 15 lines and 15 points, the Weddle quartic surface and the
  Kummer surface.
\end{abstract}

\maketitle





\section*{Introduction}

The moduli space of vector bundles of rank $2$ on a smooth projective curve $C$ has been much
studied since its construction in the sixties. it has beautiful connections to the geometry of the
$2\Theta$ linear system.  Since the Kummer variety is contained in the above moduli space (with
fixed determinant of degree 0), much of the classical Kummer geometry has vector bundle theoretical
interpretations. Surprisingly, a similar study of the $3\Theta$ linear system also reveals
connections to classical complex algebraic geometry.

Let $C$ be a smooth complex projective curve of genus~2 and $J^1$ the space of divisors of degree 1
on $C$.  The variety $J^1$ has a canonical Riemann theta divisor $\Theta$.  We denote by $\SU_C(3)$
the moduli space of semi-stable rank-3 vector bundles on $C$ with trivial determinant.  This
projective variety of dimension~$8$ is the double cover of $|3\Theta| \cong \PP^8$ branched along a
hypersurface of degree $6$ that we call the \emph{Coble sextic} and denote by $\calC_6$.  Moreover
$J^1$ embeds naturally into $|3\Theta|^* \cong \check{\PP}^8$, and there is a unique cubic
hypersurface $\calC_3$ in $\check{\PP}^8$ singular along the embedded $J^1$.  We call it the
\emph{Coble cubic}.  The main result of this paper is the ``global'' duality:

\begin{varthm}{Theorem \ref{thm:duality}}
The Coble hypersurfaces $\calC_3$ and $\calC_6$ are dual.
\end{varthm}

This was first conjectured by Dolgachev and mentioned by Laszlo in \cite{Las96}.  It was finally
proved by Ortega Ortega in her thesis \cite{Ort05}, with the use of some computer calculations.  We
give here a different ``computer-free'' proof.  The duality then allows us to deduce a non-abelian
Torelli result.

\begin{varthm}{Corollary \ref{cor:torelli}}
  Let $C$ and $C'$ be two smooth projective curves of genus~2.  If $\SU_C(3)$ is isomorphic to
  $\SU_{C'}(3)$, then $C$ is isomorphic to $C'$.
\end{varthm}

Using the duality and a geometric study of $\SU_C(3)$, we can recover a number of classical results
related to curves of genus 2, but in the context of vector bundles.

\section*{Acknowledgements}

I would like to express my sincere thanks to Igor Dolgachev for all the support, guidance and
encouragement during my research.  His great insights contributed to my fond interest in the
subject.
  
For the material covered in this treatise, I am greatly indebted to Igor Dolgachev for all the long
talks we had, but also to Alessandro Verra, S\l awek Rams, Angela Ortega Ortega, Mihnea Popa, Paul
Hacking and Ravi Vakil.


\vspace{18pt}
\section{Preliminaries: Definitions and Notations}
\label{prelim}

\vspace{12pt}
\subsection{The moduli space of vector bundles: generalities} \label{sec:moduli}

Let $E$ be an algebraic (or holomorphic) vector bundle of rank $r$ on a smooth projective curve $C$
of genus $g(C)=g$.  We define its \emph{slope of $E$} to be the number
\begin{equation*}
\mu(E) = \frac{\deg(E)}{r} = \frac{\deg(\det(E))}{r}.
\end{equation*}

A vector bundle $E$ is said to be \emph{semi-stable} (resp. \emph{stable}) if for any proper
subbundle $F$ the following inequality holds
\begin{equation*}
  \mu(F) \leq \mu(E) \quad (\text{resp. } \mu(F) < \mu(E)).
\end{equation*}

This notion of stability leads to the Geometric Invariant Theoretical construction of moduli spaces
of vector bundles, under the following equivalence relation \cite{Ses67}: every semi-stable vector
bundle $E$ admits a \emph{Jordan-H\"older filtration}
\begin{equation*}
0 = E_0 \subsetneq E_1 \subsetneq \dots \subsetneq E_{k-1} \subsetneq E_{k} = E,
\end{equation*}
such that each successive quotient $E_i/E_{i-1}$ is stable of slope equal to $\mu(E)$, for
$i=1,\ldots,k$.  We call
\begin{equation*}
\gr(E) = \bigoplus_{i=1}^k E_i/E_{i-1}
\end{equation*}
the \emph{graded bundle} associated to $E$.  Finally, two semi-stable vector bundles $E$ and $E'$ on
$C$ are said to be \emph{$S$-equivalent} if $\gr(E) \cong \gr(E')$.  We write $E \sim_{S} E'$.  In
particular, two stable bundles $E$ and $E'$ are $S$-equivalent if and only if they are isomorphic.

So we denote by $\U_C(r,d)$ the \emph{moduli of $S$-equivalence classes of semi-stable vector
  bundles on $C$ of rank $r$ and degree $d$}.  If we fix a line bundle $L$ on the curve $C$, we
denote by $\SU_C(r,L)$ the \emph{moduli space of $S$-equivalence classes of semi-stable vector
  bundles on $C$ of rank $r$ and fixed determinant $L$}.  The singular locus of these moduli spaces
corresponds exactly to decomposable bundles, i.e. \emph{strictly semi-stable bundles}, except when
$g(C)=r=2$ and $d$ is even.

Tensoring by a line bundle $M$ on $C$ induces an isomorphism
\begin{equation*}
\SU_C(r,L) \xrightarrow{\sim} \SU_C(r,L\otimes M^r),
\end{equation*}
so when we do not want to specify $L$, we will just write $\SU_C(r,d)$, where $d=\deg(L)$.

In \cite{DN89}, Drezet and Narasimhan prove that these spaces are locally factorial and
describe the Picard groups.  For the moduli space with trivial determinant
$\SU_C(r)=\SU_C(r,\OO_C)$, if we fix a general line bundle $L$ of degree $g-1$, then the set
\begin{equation*}
  \Delta_L = \{ E \in \SU_C(r) \,:\, h^0(C,E\otimes L)>0 \}
\end{equation*}
is a divisor on $\SU_C(r)$ whose isomorphism class does not depend on $L$.  We write
\begin{equation*}
  \Thgen = \OO_{\SU_C(r)}(\Delta_L)
\end{equation*}
for the corresponding line bundle (isomorphism class).  Then the Picard group of $\SU_C(r)$ is
infinite cyclic and generated by $\Thgen$.

This ample generator is called the \emph{generalized theta divisor} (or \emph{determinant bundle})
for it generalizes the traditional notion of theta divisor on Jacobians of curves.  Indeed the
variety $J^{g-1}$, seen as a space of line bundles of degree $g-1$ on $C$, has a \emph{canonical
  Riemann theta divisor} $\Theta$ defined as
\begin{equation*}
  \Theta = \{ L \in J^{g-1} \, :\, h^0(C,L)>0 \}.
\end{equation*}

For any $E \in \SU_C(r)$, we define
\begin{equation*}
D_E = \{ L \in J^{g-1} \,:\, h^0(C,E\otimes L)>0 \} \subset J^{g-1}.
\end{equation*}
It is known that $D_E$ is either the whole space $J^{g-1}$ or a divisor of the linear system
$|r\Theta|$.  The former case only happens for special $E \in \SU_C(r)$, so we get a rational map
$\Phi_r$:
\begin{equation}
  \label{eq:Phi}
  \begin{split}
    \Phi_r\ :\ \SU_C(r) &\dashrightarrow |r\Theta|,\\
    E &\longmapsto D_E=\{ L \in J^{g-1} \,:\, h^0(C,E\otimes L)>0 \}.
  \end{split} 
\end{equation}
This map is a canonical description of the map defined by the linear system $|\Thgen|$.  This
follows from a theorem of Beauville, Narasimhan and Ramanan \cite{BNR89} which states
that there is a canonical isomorphism
\begin{equation*}
H^0(\SU_C(r),\Thgen)^{*} \cong H^0(J^{g-1}(C),r\Theta).
\end{equation*}

\begin{eg} \label{eg:Phi2} For a curve $C$ of genus~2, the rational map $\Phi_2:\ \SU_C(2) \to
  |2\Theta|$ is an isomorphism (see \cite{NR69}), i.e. $\SU_C(2) \cong \PP^3$.
\end{eg}

\vspace{12pt} \subsection{Origins and motivations: The Coble quartic} \label{sec:c4}

Before we start discussing vector bundles of rank 3, we look at what is already known in rank 2 and
the main motivating result.

A concrete geometric description of moduli spaces of vector bundles with fixed determinant on a
given curve $C$ is known in the following cases:
\begin{itemize}
\item $g(C)=0$.  $\SU_C(r,d)$ is either empty or just a point (when $r|d$).
\item $g(C)=1$.  $\SU_C(r,d) \cong \PP^{h-1}$, where $h=(r,d)$ \cite{Ati57b,Tu93}.
\item $g(C)=2$, $r=2$.  $\SU_C(2)\cong\PP^3$ and $\SU_C(2,1)$ is isomorphic to the intersection of
  2 quadrics in $\PP^5$ \cite{NR69,New68}.
\end{itemize}

Then let $C$ be non-hyperelliptic of genus 3 and $r=2$.  The divisor $2\Theta$ on the Jacobian $J^2$
defines a map from $J^2$ to $|2\Theta|^*$ whose image is the Kummer variety $\mathcal{K}_2$.  In
\cite{Cob29}, Coble shows that there is a unique quartic hypersurface---now called the \emph{Coble
  quartic}---in $|2\Theta|^*\cong\PP^7$ singular exactly along $\mathcal{K}_2$.  On the other hand,
Narasimhan and Ramanan \cite{NR87} study the natural map
\begin{equation*}
  \Phi_2:\ \SU_C(2) \to |2\Theta|
\end{equation*}
and prove that it embeds $\SU_C(2)$ as a quartic hypersurface of $|2\Theta|$ singular exactly along
the Kummer variety.  So by uniqueness, it follows that $\SU_C(2)$ is isomorphic to the Coble
quartic.  Using the moduli interpretation (and the Wirtinger duality), Pauly proves that the Coble
quartic is self-dual \cite{Pau02}.

The beautiful connections found in the geometry of the Coble quartic provide the inspiration for the
present account.

\vspace{12pt} \subsection{The Coble sextic and the Coble cubic} \label{sec:c6}

Let $C$ be a smooth projective curve of genus 2, therefore hyperelliptic.  We wish to study the
moduli space $\SU_C(3)$ of rank-3 vector bundles on $C$ with trivial determinant. As usual, $\Theta$
denotes the canonical Riemann theta divisor of $J^1$ and $\Thgen$ the generalized theta divisor.
The main tool to study $\SU_C(3)$ is the map $\Phi_3$ of \eqref{eq:Phi}.

\begin{thm} \label{thm:Phi3}
The map $\Phi_3:\ \SU_C(3) \to |3\Theta| \cong \PP^8$ is a finite map of degree $2$.
\end{thm}

A first unpublished proof was given by Butler and Dolgachev using the Verlinde formula, then Laszlo
produced another beautiful proof in \cite{Las96} by making a Hilbert polynomial computation.


Let $\calI_{J^1}$ be the ideal sheaf corresponding to the tricanonical embedding $J^1
\hookrightarrow \check{\PP}^8$.  Since $J^1$ is projectively normal \cite{Koi76}, it follows that
\begin{equation} \label{eq:9quadrics}
  h^0(\check{\PP}^8,\calI_{J^1}(2)) = h^0(\check{\PP}^8,\OO(2)) - h^0(J^1,\OO_{J^1}(6\Theta)) = 9.
\end{equation}

\begin{rmk} \label{rmk:notation} 
  We use the notation $\PP^8$ for $|3\Theta|$, and therefore $\check{\PP}^8$ for for the dual
  projective space $|3\Theta|^*$.
\end{rmk}

In \cite{Cob17}, Coble produces nine quadrics that cut out the embedded Jacobian $J^1$, even
scheme-theoretically \cite{Bar95}.  More remarkably, it turns out that the nine quadrics are the
nine partial derivatives of a cubic polynomial, leading to the conclusion that there is a unique
cubic hypersurface in $\check{\PP}^8$ singular exactly along $J^1$ \cite{Bea03}.

\begin{dfn} \label{dfn:C3} The cubic hypersurface singular exactly along $J^1$ is called the
  \emph{Coble cubic}.  We will denote it by $\calC_3$.
\end{dfn}

Back to the double cover $\Phi_3$, an easy computation using the Hurwitz formula and the fact that
the canonical bundle of $\SU_C(3)$ is $(\Thgen)^{-6}$ shows that the branch locus of
$\Phi_3$ is a sextic hypersurface in $|3\Theta| \cong \PP^8$.

\begin{dfn} The branch divisor of $\Phi_3:\ \SU_C(3) \to |3\Theta|$ is called the \emph{Coble sextic}.
  We will denote by it $\calC_6$.
\end{dfn}

The name comes from Dolgachev's conjecture that this branch locus is the dual variety of the Coble
cubic, a statement clearly motivated by the Coble quartic (see Section \ref{sec:c4}) and its
self-duality.


\vspace{18pt} \section{A First Vector Bundle Interpretation of Classical Geometry}

\vspace{12pt} \subsection{Natural group actions on $|3\Theta|$ and equations of the Coble cubic}

We will now describe two natural group actions on $|3\Theta|$ (and also on $|3\Theta|^*$).  First,
there is the involution of the double cover map $\Phi_3$.  The standard adjunction involution of
$J^1$, $L \mapsto \omega_C \otimes L^{-1}$, induces involutions of the vector spaces
$H^0(J^1,\OO(3\Theta))$ and $H^0(J^1,\OO(3\Theta))^*$, and then of the projectivizations $|3\Theta|
\cong \PP^8$ and $|3\Theta|^* \cong \check{\PP}^8$.  We denote by $\tau$ both the involutions of
$|3\Theta|$ and $|3\Theta|^*$.  Let also $\tau'$ be the involution of $\SU_C(3)$ given by $E \mapsto
\tau'(E)=E^{*}$, where $E^{*}$ denotes the dual vector bundle of $E$, and let $h$ be the
hyperelliptic involution of $C$.  Then, the double cover involution $\sigma$ is (see for instance
\cite{Ort05})
\begin{equation*}
  \sigma = \tau' \circ h^* = h^* \circ \tau' :\ E \mapsto h^* E^*,
\end{equation*}
that is, the ramification locus of $\SU_C(3)$ corresponds exactly to
\begin{equation*}
  \{ E \in \SU_C(3)\,:\, \sigma(E)=h^* E^* \sim_{S} E \} \cong \calC_6.
\end{equation*}

By Riemann-Roch, we see that $\Phi_3$ is $\tau$-equivariant, i.e.
\begin{equation} \label{eq:tauphi=phitau'}
  \tau \circ \Phi_3 = \Phi_3 \circ \tau'.
\end{equation}
This implies in particular that if we
  identify the branch locus and the ramification locus $\calC_6$ of the double cover $\Phi_3$ of
  $\PP^8$, then
\begin{equation} \label{eq:intersectioninC6}
  \Fix(\tau) \cap \calC_6 \cong \Fix(\tau') \cap \calC_6,
\end{equation}
where $\Fix(\tau)$ (resp. $\Fix(\tau')$) denotes the fixed locus of $\tau$ (resp. $\tau'$).  So we
can also describe points of $\Fix(\tau)$ or $\calC_6 \subset \PP^8$ as vector bundles.

Another natural group acting on $|3\Theta| \cong \PP^8$ is the group $J_3$ of 3-torsion points of
the Jacobian $J$.  We know that $J_3$ is symplectically isomorphic to the group $(\FF_3)^4$, where
$\FF_3$ denotes the cyclic group of order 3.  The choice of such an isomorphism is called a
\emph{level-3 structure} and corresponds to the choice of a nice basis for $H^0(J^1,\OO(3\Theta))$ in
the following sense (see for instance \cite{LB92}):
\begin{thm} \label{thm:Phiisom}
  Let us fix a symplectic isomorphism $\phi:\ J_3 \to (\FF_3)^4$, where the symplectic structure on
  $J_3$ is the Weil pairing defined by the cup product on $H^1(C,\FF_3)\cong J_3$.  Then there
  exists a unique isomorphism $|3\Theta| \xrightarrow{\sim} \PP^8$ which is $\phi$-equivariant with
  respect to the action of $J_3$ on $|3\Theta|$ and that of $(\FF_3)^4$ on $\PP^8$ under the
  Schr\"{o}dinger representation.
\end{thm}
From the Schr\"{o}dinger representation of $(\FF_3)^4$, coordinates on $\CC^9$ can be written $X_b$,
for $b \in (\FF_3)^2$.  With these coordinates and using the quadrics derived in \cite{Cob17}, the Coble
cubic $\calC_3$ is defined by the polynomial
\begin{equation} \label{eq:C3}
  \begin{split}
    \frac{\alpha_0}{3} \sum_{b \in (\FF_3)^2} X_b^3 \quad
    &+ 2\alpha_1 \left( X_{00}X_{01}X_{02}+X_{10}X_{11}X_{12}+X_{20}X_{21}X_{22} \right) \\
    &+ 2\alpha_2 \left( X_{00}X_{10}X_{20}+X_{01}X_{11}X_{21}+X_{02}X_{12}X_{22} \right) \\
    &+ 2\alpha_3 \left( X_{00}X_{11}X_{22}+X_{01}X_{12}X_{20}+X_{10}X_{21}X_{02} \right) \\ 
    &+ 2\alpha_4 \left( X_{00}X_{12}X_{21}+X_{01}X_{10}X_{22}+X_{02}X_{11}X_{20} \right),
  \end{split}
\end{equation}
where $\alpha_0$, \dots, $\alpha_4$ are parameters for the genus-2 curve $C$.  Moreover, $J_3$ acts
by tensor product on $\SU_C(3)$ and it is easy to prove the following from \eqref{eq:Phi} and
\eqref{eq:C3}.
\begin{prop} \label{prop:J3-equiv} The double cover map $\Phi_3$ is $J_3$-equivariant.  In addition,
  the Coble sextic $\calC_6$ and the Coble cubic $\calC_3$ are both $J_3$-invariant.
\end{prop}

\vspace{12pt} \subsection{The Igusa-Segre quartic}
\label{sec:igusa}

The fixed locus of the involution $\tau$ of $|3\Theta|$ is a disjoint union of projective spaces:
\begin{equation} \label{eq:fixtau} 
  \Fix(\tau) = \Fix(\tau)_{+} \sqcup \Fix(\tau)_{-}, \quad
  \text{i.e } \Fix(\tau) = \PP^4_{+} \sqcup \PP^3_{-}.
\end{equation} 
In this section, we state some results about the geometry of the intersections $\mathcal{V} =
\calC_6 \cap \PP^4_{+}$ and $\mathfrak{H} = \calC_6 \cap \PP^3_{-}$.  Notice first that the two
fixed component $\PP^4_{+}$ and $\PP^3_{-}$ are not contained in $\calC_6$ (\cite[Proposition
III.1]{Ngu05} and an easy consequence of its proof).

Let us fix a level-$2$ structure on $J$ or rather $J^1$.  This is a symplectic isomorphism
\begin{equation*}
  \xi:\ J_2 \to (\FF_2)^4,
\end{equation*}
where $J_2$ denotes the group of 2-torsion points of $J$.  Moreover, $J_2$ acts on the moduli space
$\SU_C(2)$ of rank-2 bundles on $C$ by tensor product.  As in Theorem \ref{thm:Phiisom} and Example
\ref{eg:Phi2}, the level-2 structure~$\xi$ determines an equivariant isomorphism $\SU_C(2)
\xrightarrow{\sim} \PP^3$ that we embed into $\mathcal{V}$ by defining
\begin{equation*}
  V_0 = \{\OO_C \oplus F \st F \in \SU_C(2) \} \cong \PP^3 \subset \mathcal{V}. 
\end{equation*}

Another way to produce vector bundles in $\mathcal{V}$ is to take symmetric powers of vector bundles
of $\SU_C(2)$.  The image of the map
\begin{equation*}
    \Sym^2:\ \SU_C(2) \to \SU_C(3),\  F \mapsto \Sym^2 F
\end{equation*}
turns out to be to a well-known quartic hypersurface of $\PP^4$, the \emph{Igusa-Segre quartic}.
This quartic, denoted $\mathcal{I}_4$, has an interpretation as the Satake compactification of the
moduli space $\mathcal{A}^*_2(2)$ of abelian surfaces with level-$2$ structure \cite{vdG82}.

\begin{thm}[\cite{NgRa03}] \label{thm:VNR}
  The scheme $\mathcal{V}$, of degree $6$ in $\PP^4_{+}$, is the union of the Igusa-Segre quartic
  $\mathcal{I}_4$ and the double hyperplane $V_0$.  Moreover, the hyperplane $V_0$ is tangent to
  $\mathcal{I}_4$ at the point corresponding to the trivial vector bundle $\OO_C^{\oplus 3}$ from
  the $\SU_C(3)$ perspective, but also to the point $(J,\xi)$ in the moduli space $\mathcal{I}_4 =
  \mathcal{A}^*_2(2)$.
\end{thm}

The geometry of $\mathcal{I}_4$ is beautiful and well known.  Its singular locus consists of 15
lines meeting in 15 points (or nodes), all fitting in a symmetric $(15_3)$-configuration: on each
line there are exactly 3 nodes and through each node pass exactly 3 lines.  However, we can identify
the 15 lines (and 15 nodes) in terms of vector bundles.  For each $\epsilon \in J_2 - \{0\}$, we set
\begin{equation*}
V_{\epsilon} = \{ L_{\epsilon} \oplus F \,:\, F \in \SU_C(2,\epsilon) \text{ and }
t_{\epsilon}(F) = F \} \cong \PP^1 \sqcup \PP^1,
\end{equation*}
where $t_{\epsilon}$ is the translation morphism of $J^1$ acting on $\SU_C(2,\epsilon)$ by tensor
product.

So $V_{\epsilon}$ is the disjoint union of two lines: one is in $\PP^4_{+}$ and the other in
$\PP^3_{-}$.  Moreover, the 15 lines lying in $\PP^4_{+}$ are exactly the 15 lines of the
singular locus of the Igusa-Segre quartic $\mathcal{I}_4$.
This in turn allows us to understand the other intersection.

\begin{thm}[\cite{NgRa03}] \label{thm:hexahedron}
  The surface $\mathfrak{H}=\calC_6 \cap \PP^3_{-}$ is a general hexahedron, i.e. the union of 6
  planes in $\PP^3_{-}$ in general position.
\end{thm}

Interestingly, intersecting the Coble sextic with $\PP^4_{+}$ or $\PP^3_{-}$ enables us to recover
the original curve $C$.  In $\PP^4_{+}$ indeed, there is a natural Kummer surface
\begin{equation} \label{eq:K'}
  \mathcal{K'} = \{ \OO_C \oplus L \oplus L^{-1} \st L \in J \} = \calI_4 \cap V_0, 
\end{equation}
which together with the tangency point---corresponding to $(J,\xi)$ from the moduli
interpretation of $\calI_4$---completely determines $C$.  Retrieving $C$ from $\PP^3_{-}$ will follow
from Theorem \ref{thm:dualhexahedron}.


\vspace{18pt} \section{The Duality of the Coble Hypersurfaces}
\label{sec:duality}

The main result of this section is Theorem \ref{thm:duality}.  Our proof will not make use of
computer calculations, unlike A.  Ortega Ortega's \cite{Ort05}.

\vspace{12pt} \subsection{The degree of the singular locus $\Sigma$}

We first analyze the singular locus $\Sigma = \Sing(\calC_6)$ of the Coble sextic.  Since the target
space of the double cover $\Phi_3$ is smooth (just $\PP^8$), we know that the singular locus
$\Sigma$ of the branch divisor is exactly the singular locus $\Sigma'$ of the covering space
$\SU_C(3)$ corresponding to strictly semi-stable vector bundles.  We will keep the two notations,
$\Sigma$ and $\Sigma'$, in order to make clear in what space we are.

\begin{lem} \label{lem:det=P3-fibration} 
  The determinant map 
  \begin{equation*}
    \det:\ \U_C(2,0) \to J,\ F \mapsto \det(F)
  \end{equation*}
  is a \mbox{$\PP^3$-fibration}: for $a \in J$, the fiber over $a$ is $\SU_C(2,a) \cong \PP^3$.  In
  particular, $\U_C(2,0)$ is a smooth variety of dimension $5$.
\end{lem}

\begin{proof}
  The map
\begin{equation} \label{eq:psi}
  \begin{split}
    \psi:\, \SU_C(2) \times J  &\xrightarrow{16:1} \U_C(2,0) \\
    (F,L) &\longmapsto F \otimes L
  \end{split}
\end{equation}
is an \'etale covering.  Indeed, if $F\otimes L\cong F'\otimes L'$, then we take the determinants
and get $L'\otimes L^{-1} = \epsilon \in J_2$.  So $L' = L \otimes \epsilon$ and $F' = F \otimes
\epsilon$.  Therefore $\U_C(2,0)$ is the quotient of the trivial projective bundle $\SU_C(2) \times
J$ under the proper and discontinuous diagonal action of $J_2$.  Thus the conclusion follows.
\end{proof}


Since $\U_C(2,0)$ is smooth, it follows the map
\begin{equation}
  \label{eq:nu}
  \nu:\ \U_C(2,0) \to \Sigma',\ F \mapsto F \oplus \det(F)^*
\end{equation}
is a resolution of singularities of $\Sigma'$, as it is clearly surjective and injective on the open
locus of stable bundles of $\U_C(2,0)$.  Therefore via $\nu$ we get:
\begin{equation} \label{eq:degSigma1}
\deg(\Sigma) = \deg_{\PP^8}(\ \Sigma \cdot H^5 \ ) = \deg_{\U_C(2,0)}(\ \nu^*(\Thgen)^5\ )
\end{equation}
where $H$ is the class of a hyperplane in $\PP^8$.  Indeed, since $\Sigma = \Phi_3{}_* \Sigma'$, we
apply the projection formula and use the resolution map $\nu$ to see that
\begin{equation*}
  \Sigma \cdot (H^5) = \Sigma' \cdot (\Phi_3^* H)^5 = \Sigma' \cdot (\Thgen)^5 =
  \nu_*U_C(2,0) \cdot (\Thgen)^5.
\end{equation*}

But for a fixed line bundle $L \in J^1$, $\Delta_{L} = \{ E \in \SU_C(3) \st H^0(C,E \otimes L) \neq
0 \}$ is a divisor representing $\Thgen$, so
\begin{equation} \label{eq:support}
  \begin{split}
    \nu^*(\Delta_{L}) &= \{ F \in \U_C(2,0) \st H^0(C,(F \otimes L) \oplus (\det(F)^* \otimes L))
    \neq 0 \}, \\
    &= \{ F \in \U_C(2,0) \st H^0(C,F\ \otimes L) \neq 0 \} \\
    & \hspace{2cm} \cup \{ F \in \U_C(2,0) \st H^0(C,\det(F)^* \otimes L) \neq 0) \}.
  \end{split}
\end{equation}
To deal with this, we define the map $\pi$ as the following composition:
\begin{equation*}
  \xymatrix{ \U_C(2,0) \ar[r]^{\pi} \ar[d]_{\det} & J^1 \\
    J \ar[r]^{-1}_{\cong} & J \ar[u]_{\otimes L}^{\cong} }
\end{equation*}
so that $\pi$ is also a $\PP^3$-fibration.  The moduli space $\U_C(2,0)$ has a generalized
theta divisor $\Thgen_{\U}$, associated to the divisor
\begin{equation*}
  \Delta'_L = \{ F \in U_C(2,0) \st H^0(C,F \otimes L) \neq 0 \} \subset \U_C(2,0).
\end{equation*}
So we see that at the divisorial level (or set-theoretically) $\nu^*(\Delta_L) = \Delta'_L \cup
\pi^*(\Theta)$, and as an isomorphism class of line bundles on $\U_C(2,0)$,
\begin{equation} \label{eq:nuThgen}
  \nu^*(\Thgen) = \Thgen_{\U} + \pi^*(\Theta).
\end{equation}

\begin{prop} \label{prop:degSigma}
  The degree of the singular locus $\Sigma$ of $\calC_6$ in $\PP^8$ is
  \begin{equation*}
    \deg(\Sigma) = 45.
  \end{equation*}
\end{prop}

\begin{proof}
  Putting \eqref{eq:nuThgen} and \eqref{eq:degSigma1} together, we obtain
  \begin{align*}
    \deg(\Sigma) &= \deg_{\U_C(2,0)}(\ [\nu^*(\Thgen)]^5\ ) = \deg_{\U_C(2,0)}\left(\ [\Thgen_{\U}]
      + [\pi^*(\Theta)]\ \right)^5 \\
    &= \sum_{i=0}^{5} \bino{5}{i} \deg_{\U_C(2,0)}\left(\ [\Thgen_{\U}]^i \cdot
      [\pi^*(\Theta)]^{5-i}\ \right) \\
    &= \deg_{\U_C(2,0)}(\ [\Thgen_{\U}]^5 \ ) + 5 \deg_{\U_C(2,0)}(\ [\Thgen_{\U}]^4
    \cdot [\pi^*(\Theta)]\ ) \\
    &\phantom{= \deg_{\U_C(2,0)}(\ [\Thgen_{\U}]^5 \ )\ } + 10 \deg_{\U_C(2,0)}(\ [\Thgen_{\U}]^3
    \cdot [\pi^*(\Theta)]^2\ ).
  \end{align*}
  There are three terms in the sum.  The technical albeit relatively easy computations can be found
  in \cite[Lemmas V.6, V.7, V.8]{Ngu05} and make use of Lemma \ref{lem:det=P3-fibration} and the
  \'{e}tale covering $\psi$ of \eqref{eq:psi}.  We hence gather
  the terms to find that
  \begin{equation*}    
    \deg(\Sigma) = 1 \times 5 + 5 \times 4 + 10 \times 2 = 45. \qedhere
  \end{equation*}
\end{proof}

\vspace{12pt} \subsection{A map given by quadrics}

Let $G$ be the homogeneous cubic polynomial defining the Coble cubic $\calC_3$.  The motivation
here is to interpret in terms of vector bundles the dual map
\begin{align*}
  \calD:\ \calC_3 \subset |3\Theta|^* &\dashrightarrow |3\Theta| \\
  p  &\longmapsto T_p(\calC_3) = \left[ \frac{\partial G}{\partial X_i} \right]_{i=0,\dots8}
\end{align*}
given by quadrics.  So we are trying to construct a rational map $\Psi:\ \calC_3 \dasharrow
\SU_C(3)$.

We know that $C \cong \Theta \subset J^1$.  So for every $a \in J$, we write
\begin{equation*}
C_a = \Theta + a \subset J^1 \quad \text{and} \quad \vartheta_a = \OO_{C_a}(\Theta|_{C_a}).
\end{equation*}
Note that $\vartheta_a$ is a line bundle of degree 2 on $C_a$ so by Riemann-Roch, we get
$h^0(C_a,\vartheta_a{}^3) = 5$ and we denote by $\check{\PP}^4_a = |\vartheta_a{}^3|^*$ the linear
span of $C_a$ in $|3\Theta|^* \cong \check{\PP}^8$.

\begin{prop} \label{prop:PP4a in C3}
  The span $\check{\PP}^4_a$ of $C_a$ in $|3\Theta|^*$ lies in $\calC_3$.
\end{prop}

\begin{proof}
  Suppose $\check{\PP}^4_a \nsubseteq \calC_3$, then $\check{\PP}^4_a \cap \calC_3 = V_3$ is a cubic
  threefold of $\check{\PP}^4_a$.  Since $\calC_3$ is singular exactly along $J^1$, then $C_a
  \subset \Sing(V_3)$, so a secant line $\ell$ to $C_a$ must lie in $V_3$.  Therefore the secant threefold
  $\text{Sec}(C_a)$ ($C_a$ does not lie in a plane) is contained in $V_3$.  We
  will get a contradiction by showing that $\deg \text{Sec}(C_a) = 8$.  Indeed, let $\ell$ be a
  general line in $\check{\PP}^4_a$, it intersects $\text{Sec}(C_a)$ at $d$ points.  Therefore, when
  we project from $\ell$, $C_a$ is mapped to a plane sextic curve of geometric genus two with $d$
  nodes.  Since the arithmetic genus of a plane sextic curve is 10, we see that $d=8$.
\end{proof}

Let $x \in \check{\PP}^4_a - C_a$.  It corresponds to a hyperplane $V_x$ in $H^0(C_a,\vartheta_a{}^3)$:
\begin{equation*}
0 \xrightarrow{} V_x \xrightarrow{j_x} H^0(C_a,\vartheta_a{}^3) \xrightarrow{x} \CC \xrightarrow{}
0.
\end{equation*}
Since $\vartheta_a{}^3$ is very ample on $C_a$ and $x \notin C_a$, $V_x$ generates
$\vartheta_a{}^3$.  We write down the evaluation exact sequences shown in Fig. \ref{fig:commdiag},
\begin{figure}
  \begin{equation*} 
    \xymatrix{
      0 \ar[r] & E_x \ar[r] \ar[d]_i & V_x \otimes \OO_{C_a} \ar[r]^{e_x} \ar[d]_{j_x} & 
      \vartheta_a{}^3 \ar[r] \ar@{=}[d] & 0 \\
      0 \ar[r] & M \ar[r] \ar[d] & H^0(C_a,\vartheta_a{}^3) \otimes \OO_{C_a} \ar[r]^{\qquad e} 
      \ar[d]_{x} & \vartheta_a{}^3 \ar[r] & 0 \\
      & \OO_{C_a} \ar@{=}[r] & \OO_{C_a} & & }
  \end{equation*}
  \caption{Commutative diagram 1}
  \label{fig:commdiag}
\end{figure}
where $i$ is an inclusion and the lower row comes from the snake lemma.  The sheaves $E_x$ and $M$
are locally free so we see them as vector bundles of rank 3 and 4 respectively and degree $-6$, so
\begin{align*}
  \mu(E_x) &= -2, & \mu(M) &= -3/2.
\end{align*}

\begin{lem}
  The vector bundle $E_x$ is semi-stable.
\end{lem}

\begin{proof}
  Suppose $F$ is a subbundle of $E_x$.  Then it is also a subbundle of $M$, but $M$ is known to be
  stable \cite{EL92} because $\deg(\vartheta_a{}^3)=6$.  So $\mu(F) < \mu(M) = -3/2$,
  therefore $\mu(F) \leq -2$ because $F$ is of rank 1 or 2, i.e $\mu(F) \leq \mu(E_x)$.
\end{proof}

In particular, $E_x(\vartheta_a)$ is semi-stable (because $E_x$ is), of rank 3, and has trivial
determinant $\OO_{C_a}$.  It fits in the twisted evaluation sequence
\begin{equation} \label{eq:ses}
0 \to E_x(\vartheta_a) \to V_x \otimes \vartheta_a \to \vartheta_a{}^4 \to 0.
\end{equation}
We can hence define a rational map $\Psi$ from $\check{\PP}^4_a$ to $\SU_C(3)$, regular outside of
$C_a$:
\begin{equation} \label{eq:Psi}
  \begin{split}
    \Psi:\ \check{\PP}^4_a-C_a &\to \SU_C(3) \\
    x &\mapsto  \Psi(x) = E_x(\vartheta_a).
  \end{split}
\end{equation}

We will now study this map to see that it is defined by quadrics.  By Riemann-Roch, we find that
there is a non trivial section $\OO_{C_a} \to E_x(\omega_{C_a})$.  Since the two vector bundles are
of degree 0, this morphism is injective and the quotient is also a vector bundle.  So when we twist
by $\vartheta_a \otimes \omega_{C_a}^{-1}$, we obtain the following short exact sequence of vector
bundles, all of degree 0:
\begin{equation*}
  0 \to \vartheta_a \otimes \omega_{C_a}^{-1} \to E_x(\vartheta_a) \to G \to 0.
\end{equation*}
It follows that $E_x(\vartheta_a) \ \sim_S \ (\vartheta_a \otimes \omega_{C_a}^{-1}) \oplus G$,
therefore
\begin{equation*}
  \det G = \omega_{C_a} \otimes \vartheta_a^{-1} \quad \text{ and } \quad \nu(G) = E_x(\vartheta_a).
\end{equation*}
But $\SU_C(2,\omega_{C_a} \otimes \vartheta_a^{-1})$ sits naturally in $\U_C(2,0)$ as a fiber of the
determinant map (see Lemma \ref{lem:det=P3-fibration}).  Moreover, it is easy to check that
\begin{equation*}
  \nu|_{\SU_C(2,\omega_{C_a} \otimes \vartheta_a^{-1})}:\ \SU_C(2,\omega_{C_a} \otimes
  \vartheta_a^{-1}) \to \Sigma' \subset \SU_C(3)
\end{equation*}
is isomorphic onto its image and that the composition
\begin{equation} \label{eq:Dnu} \Phi_3 \circ \nu:\ \SU_C(2,\omega_{C_a} \otimes \vartheta_a^{-1})
  \to \Sigma' \to \Sigma \subset |3\Theta|=\PP^8
\end{equation}
embeds $\SU_C(2,\omega_{C_a} \otimes \vartheta_a^{-1})$ as a linear subspace of $|3\Theta|$.  We
hence write $\PP^3_a$ for $\SU_C(2,\omega_{C_a} \otimes \vartheta_a^{-1})$ and see it as a subspace
of $\Sigma$ or $\Sigma'$ interchangeably.  Thus we have just proved that the map $\Psi$ actually
lands into $\PP^3_a$.

\begin{prop} \label{prop:Psi} The rational map $\Psi:\ \check{\PP}^4_a \dashrightarrow \PP^3_a$ of
  \eqref{eq:Psi} is given by a linear system of quadrics.
\end{prop}

\begin{proof}
  The degree of the linear system defining $\Psi$ is the degree of $\Psi^*(\Thgen)$.  At the
  divisorial level, if we fix $L \in J^1$, this is just
\begin{equation*}
  \Psi^*(\Delta_L) = \{ x \in \check{\PP}^4_a \st H^0(C_a,E_x(\vartheta_a)\otimes L) \neq 0 \}
\end{equation*}
where $E_x(\vartheta_a) = \Psi(x)$.  Let us choose $L \in J^1$ so that $\vartheta_a \otimes L$ is
globally generated.  By Riemann-Roch, it is easy to see that many such $L$ exist.  If we twist the
commutative exact diagram of Fig. \ref{fig:commdiag} by $\vartheta_a \otimes L$, we get the
commutative ``long exact'' diagram of Fig. \ref{fig:commdiag2}
\begin{figure}
  \begin{equation*} 
    \xymatrix{& 0 \ar[d] & 0 \ar[d] & \\
      0 \ar[r] & H^0(E_x(\vartheta_a)\otimes L) \ar[r] \ar[d] & V_x \otimes H^0(\vartheta_a \otimes L) 
      \ar[r]^{\quad e_x} \ar[d] & H^0(\vartheta_a{}^4 \otimes L) \ar@{=}[d] \\
      0 \ar[r] & H^0(M(\vartheta_a)\otimes L) \ar[r] \ar[d] & H^0(\vartheta_a{}^3) \otimes 
      H^0(\vartheta_a \otimes L) \ar[d]_{g(x)} \ar[r]^{\qquad \quad e} & H^0(\vartheta_a{}^4 \otimes L) \\
      & H^0(\vartheta_a \otimes L) \ar@{=}[r] & H^0(\vartheta_a{} \otimes L)  & }
  \end{equation*}
  \caption{Commutative diagram 2} \label{fig:commdiag2}
\end{figure}
where the cohomology groups are taken over $C_a$ and where the map $g(x)$ can be described as
follows.  The contraction
\begin{equation*}
  H^0(\vartheta_a{}^3) \otimes H^0(\vartheta_a{}^3)^* \otimes H^0(\vartheta_a \otimes L) \to
  H^0(\vartheta_a \otimes L)
\end{equation*}
defines a linear map
\begin{equation*}
  g:\ H^0(\vartheta_a{}^3)^* \to \Hom (H^0(\vartheta_a{}^3) \otimes H^0(\vartheta_a \otimes L),
  H^0(\vartheta_a \otimes L)).
\end{equation*}
We will abuse notation and write $x$ for both an element of $\PP H^0(\vartheta_a{}^3)^*$ and any of
its representatives in $H^0(\vartheta_a{}^3)^*$.  So for all $x \in H^0(\vartheta_a{}^3)^*$,
\begin{gather*}
  g(x):\ H^0(\vartheta_a{}^3) \otimes H^0(\vartheta_a \otimes L) \to H^0(\vartheta_a \otimes L),
  \quad \text{and} \\
  \ker(g(x)) = V_x \otimes H^0(\vartheta_a \otimes L).
\end{gather*}
Since $\dim(V_x \otimes H^0(\vartheta_a \otimes L)) = \dim(H^0(\vartheta_a{}^4 \otimes L)) = 8$, we
have exactly
\begin{equation*}
  \Psi^*(\Delta_L) = \{ x \in \check{\PP}^4_a \st \text{the map $e_x$ degenerates} \}.
\end{equation*}
Moreover $h^0(\vartheta_a \otimes L)=2$ and $\vartheta_a \otimes L$ is globally generated, so by the
base point free pencil trick, 
$\dim \ker(e) = h^0(\vartheta_a{}^2 \otimes L^*) = 2$.  Therefore by restricting $g$, we get
\begin{equation*}
  g':\ H^0(\vartheta_a{}^3)^* \to \Hom(H^0(\vartheta_a{}^2 \otimes L^*),H^0(\vartheta_a \otimes L))
  \cong \Hom(\CC^2,\CC^2),
\end{equation*}
and rewrite the commutative exact diagram of Fig. \ref{fig:commdiag2} as that of Fig.
\ref{fig:commdiag3}.
\begin{figure}
  \begin{equation*}
    \xymatrix{& 0 \ar[d] & 0 \ar[d] & \\
      0 \ar[r] & H^0(E_x(\vartheta_a)\otimes L) \ar[r] \ar[d] & V_x \otimes H^0(\vartheta_a \otimes L)
      \ar[r]^{e_x} \ar[d] & H^0(\vartheta_a{}^4 \otimes L) \ar@{=}[d] \\
      0 \ar[r] & H^0(\vartheta_a{}^2 \otimes L^*) \ar[r] \ar[d]_{g'(x)} & H^0(\vartheta_a{}^3) 
      \otimes H^0(\vartheta_a \otimes L) \ar[d]_{g(x)} \ar[r]^{\qquad \quad e} &
      H^0(\vartheta_a{}^4 \otimes L) \\
      & H^0(\vartheta_a{} \otimes L) \ar@{=}[r] & H^0(\vartheta_a{} \otimes L). & }
  \end{equation*}
  \caption{Commutative diagram 3} \label{fig:commdiag3}
\end{figure}
It is then clear that $e_x$ degenerates, i.e. $H^0(E_x(\vartheta_a)\otimes L) \neq 0$, exactly when
$g'(x)$ degenerates, since $h^0(\vartheta_a{}^2 \otimes L^*) = h^0(\vartheta_a \otimes L) = 2$.  So
\begin{equation*}
  \Psi^*(\Delta_L) = \{ x \in \PP H^0(\vartheta_a{}^3)^* \st \det(g'(x)) = 0 \},
\end{equation*}
is a quadric as it is the pull-back under $g'$ of the discriminant locus of
\begin{equation*}
  \Hom(H^0(\vartheta_a{}^2 \otimes L^*),H^0(\vartheta_a \otimes L)). \qedhere
\end{equation*}
\end{proof}

\vspace{12pt} \subsection{Restriction of the dual map}

Now that we have proved that the map $\Psi$ is given by quadrics, we will show that it is actually
the restriction of the dual map $\calD:\ \calC_3 \dasharrow |3\Theta|$. 

\begin{dfn}
  The dual variety of $\calC_3$, denoted $\check{\calC}_3$, is the image of the dual map $\calD$.
  We also call the dual map the rational map on $\check{\PP}^8$.  By abuse of notation, $\calD$
  denotes both the map $\calC_3 \dasharrow \PP^8$ and the map $\check{\PP}^8 \dasharrow
  \PP^8$.
\end{dfn}

By definition, the dual map or polar map $\calD$ is given by 9 quadrics passing through the singular
locus of the Coble cubic $\calC_3$.  The dimension count \eqref{eq:9quadrics} shows that $\calD$ is
given by the complete linear system $|\calI_{J^1}(2)|$.  So when restricted to $\check{\PP}^4_a$,
for $a \in J$, $\calD$ is given by quadrics in $\check{\PP}^4_a$ containing $C_a=\Theta + a \subset
J^1$.

\begin{prop}
\label{prop:restriction=c.l.series}
The restriction of the dual map $\calD|_{\check{\PP}^4_a} :\ \check{\PP}^4_a \dasharrow |3\Theta|$
is given by the complete linear system $|\calI_{C_a}(2)| \cong \PP^3$ of quadrics in
$\check{\PP}^4_a$ containing $C_a$.
\end{prop}

\begin{proof}
  Twisting a usual sequence by $\OO_{\check{\PP}^4_a}(2)$, we get the exact sequence
  \begin{equation} \label{eq:seqCa} 
    0 \to \calI_{C_a}(2) \to \OO_{\check{\PP}^4_a}(2) \to \OO_{C_a}(\vartheta_a{}^6) \to 0.
  \end{equation}
  Since $C_a$ is of genus~2, $\calI_{C_a}$ is 3-regular \cite[Theorem 2.1]{GLP83}, so
  $h^0(\check{\PP}^4_a,\calI_{C_a}(2)) = 4$ by \eqref{eq:seqCa}.  Let $\alpha:\ |\calI_{J^1}(2)| \to
  |\calI_{C_a}(2)|$ be the natural restriction map.  The proposition is equivalent to $\alpha$ being
  surjective. By contradiction, let us assume that the rank of $\alpha$ is less than 4.  So we can
  choose a basis of the 9-dimensional space $H^0(\check{\PP}^8,\calI_{J^1}(2))$ that consists of at
  least 6 quadrics that contain $\check{\PP}^4_a$ and at most 3 quadrics that do not.  Let
  $\mathcal{B}$ be the base locus restricted to $\check{\PP}^4_a$ of the latter 3 quadrics: $\dim
  \mathcal{B} \geq 1$.  The base locus of $H^0(\check{\PP}^8,\calI_{J^1}(2))$, by definition exactly
  $J^1$, must then contain $\mathcal{B}$.  But $J^1 \cap \check{\PP}^4_a = C_a$ is a curve of degree
  $6$ in $\check{\PP}^4_{a}$.  However, $C_a \supseteq \mathcal{B}$ only if $\mathcal{B}$ has
  dimension 1, in which case $\deg(\mathcal{B})=8$: a contradiction.
\end{proof}



\begin{prop} \label{prop:Psi=D'}
  The rational maps $\Psi$ and $\calD|_{\check{\PP}^4_a}$ are equal.
\end{prop}

\begin{proof}
  By definition, $\Psi$ is regular outside of $C_a$, so its base locus $B$ is contained in $C_a$.
  For the reverse inclusion, we want to show that $g'(x)$ never has maximal rank for all $x \in
  C_a$.  Elements of $H^0(C_a,\vartheta_a{}^2 \otimes L^*)$ are linear combinations of tensors $ s
  \otimes \sigma \in H^0(C_a,\vartheta_a{}^3) \otimes H^0(C_a,\vartheta_a \otimes L)$ such that
  $s\cdot\sigma \in H^0(C_a,\vartheta_a{}^4 \otimes L)$ is zero (see Fig. \ref{fig:commdiag3}).
  But on simple tensors, $g'(x)$ acts by
  \begin{equation*}
      g'(x):\ H^0(C_a,\vartheta_a{}^2 \otimes L^*) \to H^0(C_a,\vartheta_a \otimes L),\ 
      s \otimes \sigma \mapsto s(x) \cdot \sigma,
  \end{equation*}
  therefore the image of $g'(x)$ is the subspace of sections in $H^0(C_a,\vartheta_a\otimes{}L)$
  vanishing at~$x$.  However, $\vartheta_a\otimes{}L$ was assumed to be globally generated (see
  proof of Proposition \ref{prop:Psi}), so the image of $g'(x)$ is a proper subspace, i.e. $g'(x)$
  degenerates.  Thus $B=C_a$.
  So $\Psi$ is given by a linear subseries of $|\calI_{C_a}(2)|$ that has to be of dimension~3 by
  Proposition \ref{prop:Psi}.  But $\calD|_{\check{\PP}^4_a}$ is given by the complete linear series
  $|\calI_{C_a}(2)|$ which also has dimension~3.  So they are equal.
\end{proof}

One natural question about the construction of $\Psi$ is: what if we see $x \in \check{\PP}^4_a$ as
an element of $\check{\PP}^4_b$, another $\PP^4$? Since $\Psi$ (or rather $\Psi_a$) coincides with
$\calD|_{\PP^4_{a}}$, it does not matter which $\PP^4$ we choose to define $\Psi(x)$ and we can
extend the definition of $\Psi$ to the variety $\mathcal{W} = \bigcup_{a \in J} \check{\PP}^4_{a}$.
We denote the extended map $\overline{\Psi}$.

\begin{prop} \label{prop:global} 
  We have a ``global'' equality of dominant rational maps:
  \begin{equation*}
    \overline{\Psi} = \calD|_{\mathcal{W}}:\ \mathcal{W} \dashrightarrow \Sigma = \Sing(\calC_6).
  \end{equation*}
  Furthermore, $\Sigma =\Sing(\calC_6) \subset \Sing(\check{\calC}_3)$.
\end{prop}

\begin{proof}
  The equality holds as we already know.  So what is the image of $\overline{\Psi}$?  First, the
  image of $\calD|_{\check{\PP}^4_a}$ is the projective space $|\calI_{C_a}(2)|^*$.  Indeed a
  general point of $|\calI_{C_a}(2)|^*$ corresponds to a net of quadrics in $\check{\PP}^4_a$ whose
  base locus is a degree-8 curve which can therefore contain $C_a$ (of degree~6).  The residual
  curve, i.e. the general fiber of $\calD|_{\check{\PP}^4_a}$, is hence a conic.  So from the point
  of view of $\Psi$, the image of $\overline{\Psi}$ is
  \begin{equation*}
    \bigcup_{a \in J} \PP^3_a = \Phi_3 \circ \nu(\U_C(2,0)) = \Sigma.
  \end{equation*}

  Furthermore, $\mathcal{W}$ is a subvariety of $\calC_3$ by Proposition \ref{prop:PP4a in C3} which
  maps onto $\Sigma$ generically forcing one-dimensional conic fibers.  By a standard property of
  dual varieties, if some subvariety of $\calC_3$ not contained in the singular locus $J^1(C)$ is
  mapped onto a lower dimensional subvariety, then that image lies in the indeterminacy locus of the
  inverse dual map, i.e. in the singular locus of the dual variety $\check{\calC}_3$.
\end{proof}

\vspace{12pt} \subsection{Finishing the proof of the duality and non-abelian Torelli}

\begin{thm} \label{thm:duality}
  The Coble hypersurfaces $\calC_3$ and $\calC_6$ are dual.
\end{thm}

First, we compute the degree of $\check{\calC}_3$ and state a general auxiliary lemma.

\begin{prop} \label{prop:degreedual} 
  The dual variety $\check{\calC}_3$ is a hypersurface of degree 6 in $\PP^8$.
\end{prop}

\begin{proof}
  A simple Chern class computation shows that the total degree of the dual map $\calD$ is 6
  \cite{Ort05}.  So if we write $d$ for the generic degree of $\calD$, we have
  \mbox{$d \cdot \deg(\check{\calC}_3) = 6$.}
  To conclude that the degree of the variety has to be 6, we use results of Section \ref{geomfix}
  that, critically, do not rely on the duality we want to prove.  We know that the intersection of
  $\calC_3$ with $\check{\PP}^4_{+}$ \eqref{eq:fix tau} is the Segre cubic $\calS_3$ (see Section
  \ref{sec:segre-cubic}).  By the commutativity of \eqref{eq:dualscommutewithembeddings2}, the
  proper intersection $\check{\calC}_3 \cap \PP^4_{+}$ contains the dual variety of $\calS_3$, which
  is known to be the Igusa-Segre quartic $\calI_4$.  Therefore the degree of $\check{\calC}_3$ is
  greater than or equal to 4.  So the only possibility is 6.
\end{proof}

\begin{lem} \label{lem:k geq d}
  Let $\mathcal{G}$ be a group acting on the projective space $\PP^n$.  Let $V_d$ be a
  $\mathcal{G}$-invariant hypersurface of degree $d$.  Let $W_k$ be a hypersurface of degree $k$
  such that the scheme-theoretic intersection $Y=V_d \cap W_k$ is $\mathcal{G}$-invariant.  If
  $k<d$, then $W_k$ is $\mathcal{G}$-invariant.
\end{lem}

\begin{proof}
  Let $F_d$ (resp. $G_k$) be a homogeneous polynomial defining $V_d$ (resp. $W_k$).  Then the
  homogeneous ideal of $Y$ is generated by $F_d$ and $G_k$.  This ideal is $\mathcal{G}$-invariant,
  by this we mean that each homogeneous part is an invariant subspace of the vector space of
  homogeneous polynomials of fixed degree.  If $k<d$, then $G_k$ is the only form of degree $k$, so
  it has to be $\mathcal{G}$-invariant.
\end{proof}

\begin{proof}[Proof of Theorem \ref{thm:duality}]
  Let us assume that $\calC_6$ and $\check{\calC}_3$ are different.  We write \makebox[1.2\width]{$Y
    = \calC_6 \cap \check{\calC}_3$.}  Since $\calC_6$ has a singular locus of codimension~2, it is
  irreducible.  Similarly $\calC_3$ is irreducible and so is its dual variety $\check{\calC}_3$
  \cite[Proposition 1.3]{GKZ94}.  Therefore, $Y$ is connected \cite[Proposition 1]{FH79}.  We must
  then consider two cases.

  \noindent
  \textbf{Case 1:} $Y$ has a reduced component.  If we intersect the objects with a general $\PP^3$,
  we see that the surface $S=\calC_6\cap\PP^3$ has 45 $A_1$-singularities, according to the local
  analytic description of $\Sigma' \subseteq \SU_C(3)$ in \cite{Las96} and the fact that
  $\deg(\Sigma)=45$.  Since $T=\check{\calC}_3\cap\PP^3$ is different from $S$, we denote by $D$ the
  Cartier divisor of $S$ defined by the complete intersection with $T$.  It is a connected curve
  with a reduced component.  We resolve the 45 rational double points and get 45 exceptional
  $(-2)$-curves $E_1, \dots, E_{45}$ in the resolution $\pi:\ \tilde{S} \to S$.
  Now let $H$ be the pullback of the hyperplane section of $S$.  Its self-intersection in
  $\tilde{S}$ is $\deg_{\tilde{S}}(H^2) = 6$.  Since $T$ is also singular at those 45 points by
  Proposition \ref{prop:global}, the proper transform $\tilde{D}$ of $D$ under the blowup map $\pi$
  is linearly equivalent to
  \begin{equation*}
    \tilde{D} = 6H - \sum_{i=1}^{45} a_i E_i, \quad a_i \geq 2.
  \end{equation*}
  But $\pi$ is a crepant resolution, so
  \begin{equation*}
    \omega_{\tilde{S}} = \pi^*(\omega_{\PP^3} \otimes \OO_{\PP^3}(S) \otimes \OO_S) = \OO(2H).
  \end{equation*}
  We can also compute the arithmetic genus $p_a(\tilde{D})$ of $\tilde{D}$, knowing that $\tilde{D}$
  is reduced:
  \begin{equation*}
    2p_a(\tilde{D}) - 2  = \deg_{\tilde{S}}(\ (K_{\tilde{S}} + \tilde{D}) \cdot \tilde{D}\ ) =
    288 + \sum_{i=1}^{45} (-2)a_i^2 \leq 288 - 360 = -72,
  \end{equation*}
  because $a_i \geq 2$.  Therefore we see that $p_a(\tilde{D}) \leq -35$, which is not possible,
  since $\tilde{D}$ is an effective reduced connected Cartier divisor on a nonsingular surface.  So
  $S=T$.  Moreover, since the intersecting $\PP^3$ was general, $\calC_6 = \check{\calC}_3$.

  \noindent \textbf{Case 2:} $Y$ has no reduced components.  We write the decomposition of $Y$ into
  irreducible components:
  \begin{equation*}
    Y = a_1 Y_1 + a_2 Y_2 + \dots + a_m Y_m,
  \end{equation*}
  where $a_i \geq 2$ and $Y_i$ are prime Cartier divisors on $\calC_6$ of respective degrees $d_i$.
  But $Y$ is of degree 36, therefore $a_1 d_1 + \dots + a_m d_m=36$.  Since $\dim \calC_6 \geq 3$,
  we know by a Lefschetz-type Theorem (\cite[Expos\'{e} XII, Corollary 3.7]{SGA2}) that the
  restriction map $\Pic \PP^8 \xrightarrow{\sim} \Pic \calC_6$ is an isomorphism.  So each prime
  divisor $Y_i$ is cut out by a hypersurface in $\PP^8$ and it follows that 6 divides $d_i$.  Thus
  the only possible cases are
  \begin{itemize}
  \item $m=1$: $(a_1,d_1)=(2,18)$ or $(3,12)$ or $(6,6)$.
  \item $m=2$: $\{(a_1,d_1),(a_2,d_2)\} = \{(2,12),(2,6)\}$ or $\{(3,6),(3,6)\}$.
  \item $m=3$: $\{(a_1,d_1),(a_2,d_2)\,(a_3,d_3)\} = \{(2,6),(2,6),(2,6)\}$.
  \end{itemize}
  In every case, we can see that the $a_i$ have a common divisor that is either 2, 3 or 6.  So we
  can rewrite
  \begin{equation*}
    Y = 2Z \text{ or } 3Z \text{ or } 6Z.
  \end{equation*}
  A contradiction will arise from the $J_3$-invariance of $Y$.  We already know that $\calC_6$ and
  $\calC_3$ are $J_3$-invariant (Proposition \ref{prop:J3-equiv}).  Moreover, it is easy to see that
  $\calD$ is $J_3$-equivariant directly from the partial derivatives of \eqref{eq:C3}.  Therefore
  $\check{\calC}_3$ is $J_3$-invariant, and so is $Y$ and then $Z$.  If $Y=3Z$ or $6Z$, then $Z$ is
  a quadric section or a hyperplane section, but there are no $J_3$-invariant quadrics or
  hyperplanes in $\check{\PP}^8$.  Thus by Lemma \ref{lem:k geq d}, we get a contradiction.  So we
  are left with the case $Y=2Z$, and $Z$ is cut out by a $J_3$-invariant cubic.  Again, by
  Proposition \ref{prop:global}, we know that $\Sigma \subset Z$.  We intersect with the fixed
  $\PP^4_{+}$ \eqref{eq:fixtau}, so $Z \cap \PP^4_{+}$ must contain $\Sigma \cap \PP^4_{+}$.  But as
  a consequence of Theorem \ref{thm:VNR}, we know that \cite[Eq. (5.2)]{NgRa03}
  \begin{equation*}
    \Sigma \cap \PP^4_{+} = 2H \cup \{ (15_3)\text{-configuration of lines and points} \},
  \end{equation*}
  where $H$ is a hyperplane of $\PP^4_{+}$.  The configuration of lines and points does not lie in a
  hyperplane, so $\Sigma \cap \PP^4_{+}$ cannot lie in a cubic hypersurface of $\PP^4_{+}$, which
  shows that $\PP^4_{+}$ must lie in the $J_3$-invariant cubic hypersurface cutting out $Z$.
  However, by Proposition \ref{prop:P4 not subset of C3}, $\PP^4_{+}$ cannot be contained in a
  $J_3$-invariant cubic.  We hence ruled out all the possibilities.
\end{proof}

The established duality will be used in the next section to recover some classical dualities and
reinterpret them in terms of vector bundles.  But an easy corollary is the following non-abelian
Torelli theorem:

\begin{cor}[Non-abelian Torelli Theorem] \label{cor:torelli}
  Let $C$ and $C'$ be two smooth projective curves of genus two.  If $\SU_C(3)$ is isomorphic to
  $\SU_{C'}(3)$, then $C$ is isomorphic to $C'$.
\end{cor}

\begin{proof}
  Starting from $\SU_C(3)$, there is a canonical way to retrieve $C$.  We first take the ample
  generator $\Thgen$ of $\Pic(\SU_C(3))$, look at the map associated to the line bundle.  The branch
  locus of the 2-1 map has dual variety a cubic hypersurface in $\PP^8$ singular exactly along the
  principally polarized Jacobian $(J^1(C),\Theta)$, which determines $C$ by the usual Torelli
  theorem.
\end{proof}

\begin{rmk}
  The non-abelian Torelli question has been raised ever since the construction of the moduli spaces
  $\SU_C(r,d)$ ($g \geq 2$).  In \cite{MN68}, Mumford and Newstead prove the theorem in the case of
  rank 2, odd degree determinant ($d=1$), and $g\ge 2$.  It is further generalized in \cite{NR75}
  and \cite{Tyu74} to all non trivial smooth cases ($(r,d)=1$).  Then Balaji proves the theorem for
  $r=2$, $d=0$ on a curve of genus $g\ge 3$ \cite{Bal90}, before Kouvidakis and Pantev extend the
  result to any $r$ and $d$, still for $g\ge 3$ \cite{KP95}.
  
  In \cite{HR04}, Hwang and Ramanan introduce a stronger non-abelian Torelli result.  If we denote
  by $\SU_C(r,d)^s$ the moduli space of stable vector bundles, which is open in $\SU_C(r,d)$, then
  $\SU_C(r,d)^s \cong \SU_{C'}(r,d)^s$ implies that $C\cong C'$, and this for any $r$ and $d$, but
  for $g\ge 4$.  Our version of non-abelian Torelli is a new case and it can also be shown to be
  strong by adapting the same argument.
\end{rmk}


\vspace{18pt} \section{The Geometry of the Fixed Loci}
\label{geomfix}

In this section, we focus our attention to the dual maps:
\begin{equation*}
  \xymatrix{
    **[l] \calC_6 \subset \PP^8 \ar@{-->}@/^1ex/[rr]^{\calD'} && **[r] \check{\PP}^8 \supset \calC_3,
    \ar@{-->}@/^1ex/[ll]^{\calD} 
  }
\end{equation*}
especially, these dual maps restricted to the fixed loci $\Fix(\tau)$ defined in \eqref{eq:fixtau}:
\begin{equation} \label{eq:fix tau}
  \Fix(\tau) = \PP^4_{+} \sqcup \PP^3_{-} \subset |3\Theta| = \PP^8, \quad \text{ and } \quad
  \Fix(\tau) = \check{\PP}^4_{+} \sqcup \check{\PP}^3_{-} \subset |3\Theta|^* = \check{\PP}^8.
\end{equation}

\begin{notation}
  We use the following convention, based on the domain space.  For all dual maps (i.e. $\calD$,
  $\calD'$, and their restrictions), the apostrophe means that the domain is or lies in
  $\PP^8\cong|3\Theta|$ and not $\check{\PP}^8$.  The use of the subscripts \mbox{$_+$ or $_-$}
  determines to which fix locus we restrict the maps.  For instance, $d'_{-}$ will denote the
  restriction $\calD'|_{\PP^3_{-}}:\ \PP^3_{-} \dasharrow \check{\PP}^8$.  Also, maps that carry
  $\check{\quad}$ are maps that stay in $\check{\PP}^8$.
\end{notation}

\vspace{12pt} \subsection{$\tau$-equivariance and target spaces}
\label{sec:det target spaces}

In order to study the dual maps $\calD$ and $\calD'$, we will verify that they are
$\tau$-equivariant.  Although it is not the most intrinsic way to proceed, we will use an explicit
basis for the vector spaces $H^0(J^1,\OO(3\Theta))$ and its dual.  This amounts to choosing a
level-3 structure on the Jacobian $J$, just as in Theorem \ref{thm:Phiisom}.  The basis is $\{ e_b
\}_{b \in (\FF_3)^2}$ and the corresponding coordinate system $\{ X_b \}_{b \in (\FF_3)^2}$.

\begin{prop} \label{prop:Cobletauinvariant}
  \begin{enumerate}[(i)]
  \item The Coble hypersurfaces $\calC_3$ and $\calC_6$ are $\tau$-invariants.  Moreover, the
    polynomials defining $\calC_3$ and $\calC_6$ are $\tau$-invariant.
    
  \item The dual maps $\calD$ and $\calD'$ are $\tau$-equivariant.  So the fixed loci are mapped to
    the fixed loci by the dual maps.
  \end{enumerate}
\end{prop}

\begin{proof}
  (i) For $\calC_3$, this is clear from the equation \eqref{eq:C3} and the fact that $\tau \cdot X_b
  = X_{-b}$.  For $\calC_6$, let $E \in \SU_C(3)$, so $\tau'(E)=E^*$.  If $E \in \calC_6$
  (i.e. $h^*E^*=E$), then $h^*(\tau'E)^*= h^*E=E^*=\tau'(E)$, which means $\tau'(E) \in \calC_6$.
  So $\calC_6$ is $\tau'$-invariant as a ramification locus, i.e.  $\tau$-invariant as a branch
  locus, and the polynomial $F_6$ defining $\calC_6$ is either $\tau$-invariant or
  $\tau$-anti-invariant. But if $F_6$ is $\tau$-anti-invariant, then it must contain $\PP^4_{+}$ for
  if $V_+$ is the vector space such that $\PP^4_{+} = \PP(V_+)$ and if we take $a\in V_{+}$, then
  \begin{equation*}
    F_6(a) = F_6(\tau \cdot a) \stackrel{\text{def}}{=} \tau \cdot F_6(a) = -F_6(a).
  \end{equation*}
  So we see that $F_6(a)$ must be zero.  But we know that $\PP^4_{+}$ is not contained in $\calC_6$
  \cite[Proposition III.1]{Ngu05}.  Thus $F_6$ is $\tau$-invariant.

  \noindent (ii) An easy exercise left to the reader.
\end{proof}

We now describe the fixed spaces: $\PP^3_{-} \subset \PP^8=|3\Theta|$ is defined by the equations
\begin{equation*}
  X_{00} = 0 , \ X_{01}+X_{02} = 0,\ X_{10}+X_{20} = 0, \ X_{11}+X_{22} = 0,\ X_{12}+X_{21} = 0,
\end{equation*}
while $\PP^4_{+}$ is defined by
\begin{equation*}
  X_{01}-X_{02} = 0,\ X_{10}-X_{20} = 0, \ X_{11}-X_{22} = 0,\ X_{12}-X_{21} = 0.
\end{equation*}
Moreover, the inclusions of $\PP^3_{-}$ and $\PP^4_{+}$ in $\PP^8$ are given by the linear
injections $\gamma_{-}:\ \CC^4 \to \CC^9$ and $\gamma_{+}:\ \CC^5 \to \CC^9$ where
\begin{equation} \label{eq:embeddings}
  \begin{split}
    \gamma_{-}(Z_0,Z_1,Z_2,Z_3) &= (X_{00},X_{01},X_{02},X_{10},X_{11},X_{12},X_{20},X_{21},X_{22}), \\
    &= (0,Z_0,-Z_0,Z_1,Z_2,Z_3,-Z_1,-Z_3,-Z_2), \quad \text{ and} \\
    \gamma_{+}(Y_0,Y_1,Y_2,Y_3,Y_4) &= (X_{00},X_{01},X_{02},X_{10},X_{11},X_{12},X_{20},X_{21},X_{22}), \\
    &= (Y_0,Y_1,Y_1,Y_2,Y_3,Y_4,Y_2,Y_4,Y_3).
  \end{split}
\end{equation}

\begin{prop} \label{prop:map to which fixed loci} 
  \begin{enumerate}[(i)]
  \item The dual map $\calD:\ \check{\PP}^8 \dashrightarrow \PP^8=|3\Theta|$ maps the fixed loci in
    the following way:
    \begin{align*}
      \calD(\check{\PP}^4_{+}) &\subset \PP^4_{+}, & \calD(\check{\PP}^3_{-}) &\subset \PP^4_{+}.
    \end{align*}
  
  \item The dual map $\calD':\ \PP^8 \dashrightarrow \check{\PP}^8=|3\Theta|^*$ maps the fixed loci in
    the following way:
    \begin{align*}
      \calD'(\PP^4_{+}) &\subset \check{\PP}^4_{+}, & \calD'(\PP^3_{-}) &\subset \check{\PP}^3_{-}.
    \end{align*}
  \end{enumerate}
\end{prop}

\begin{proof}
  (i) Let $X_{00}, X_{01}, \dots, X_{22}$ be the homogeneous coordinates of $\PP^8=|3\Theta|$ and
  $F$ be a homogeneous $\tau$-invariant (Proposition \ref{prop:Cobletauinvariant}) polynomial of
  degree~6 defining the Coble sextic.  For notational convenience, let us write for $b \in
  (\FF_3)^2$
  \begin{equation} \label{eq:F'}
    F'_{b}(Z_0,Z_1,Z_2,Z_3) 
    = \frac{\partial F}{\partial X_{b}}(0,Z_0,-Z_0,Z_1,Z_2,Z_3,-Z_1,-Z_3,-Z_2),
  \end{equation}
  a homogeneous polynomial of degree 5. So
  \begin{equation}  \label{eq:F'_{01}=-F'_{02}}
    F'_{b}(Z_0,Z_1,Z_2,Z_3) = F'_{-b}(-Z_0,-Z_1,-Z_2,-Z_3) = - F'_{-b}(Z_0,Z_1,Z_2,Z_3),
  \end{equation}
  which proves that
  \begin{equation} \label{eq:F'_{b}=-F'_{-b}} 
    F'_{00} = 0, \quad F'_{10} = -F'_{20}, \quad F'_{11} = -F'_{22}, \quad \text{and } F'_{12} = -F'_{21}.
  \end{equation}
  Hence, the dual map restricted to $\PP^3_{-}$ is given as follows:
  \begin{align*}
    \calD' \circ \gamma_{-}:\ \PP^3_{-} &\dashrightarrow  \check{\PP}^8 \\
    [Z_0:Z_1:Z_2:Z_3] &\longmapsto [0:F'_{01}:-F'_{01}:F'_{10}:F'_{11}:F'_{12}:-F'_{10}:-F'_{12}:-F'_{11}]
  \end{align*}
  and the image clearly lies in $\check{\PP}^3_{-}$.  Then similarly, for $\PP^4_{+}$, if we write
  \begin{equation*}
    G'_{b}(Y_0,Y_1,Y_2,Y_3,Y_4) = \frac{\partial F}{\partial X_{b}}(Y_0,Y_1,Y_1,Y_2,Y_3,Y_4,Y_2,Y_4,Y_3),
  \end{equation*}
   it is easy to see that the dual map $\calD'$ sends $\PP^4_{+}$ to $\check{\PP}^4_{+}$ since
  \begin{equation*}
    G'_{01} = G'_{02}, \quad G'_{10} = G'_{20}, \quad G'_{11} = G'_{22}, \quad G'_{12} = G'_{21}.
  \end{equation*}

  \noindent (ii) The proof is completely similar, except for the fact that the dual map $\calD$ is
  given by quadrics instead of quintics, so it sends $\check{\PP}^3_{-}$ to $\PP^4_{+}$.
\end{proof}

\vspace{12pt} \subsection{The Hexahedron and the original curve}

In this section, we will study the image of the hexahedron $\mathfrak{H}$ of Theorem
\ref{thm:hexahedron} under the dual map $\calD'$ and prove that it also determines our original
curve $C$ of genus~2.

\begin{thm} \label{thm:dualhexahedron} 
  The 6 planes making $\mathfrak{H}$ correspond to the 6 Weierstra\ss{} points of our original
  curve~$C$.  That is, the six planes correspond to 6 points in the dual projective space
  $\check{\PP}^3_{-}$ of $\PP^3_{-}$ through which passes a unique twisted cubic.  On this rational
  curve, the points are projectively equivalent to the 6 Weierstra\ss\ points of our given curve
  $C$.
\end{thm}

First, let us introduce more notation.  Again, we use $Z_0,Z_1,Z_2,Z_3$ as the homogeneous
coordinates of $\PP^3_{-}$ and $F$ as a polynomial defining $\calC_6$ in $\PP^8$.  We keep the
notations of \eqref{eq:embeddings} and \eqref{eq:F'}.  Let
\begin{equation*}
  H(Z_0,Z_1,Z_2,Z_3) = F \circ \gamma_{-} = F(0,Z_0,-Z_0,Z_1,Z_2,Z_3,-Z_1,-Z_3,-Z_2)
\end{equation*}
be the equation defining the hexahedron.  Let $\delta'_{-}$ be the ``restricted'' dual rational map
\begin{equation}
  \begin{split}
    \delta'_{-}:\ \PP^3_{-}  & \dashrightarrow \check{\PP}^3_{-} \\
    [Z_0:\dots:Z_3] &\mapsto \left[ \frac{\partial H}{\partial Z_0} : \dots : \frac{\partial
        H}{\partial Z_3} \right].
  \end{split}
\end{equation}

\begin{prop} \label{prop:dualscommutewithembeddings} 
  Let $\check{\gamma}_{-}$ be the embedding of $\check{\PP}^3_{-}$ in $\check{\PP}^8$.  Then the
  following diagram commutes.
  \begin{equation*}
    \xymatrix{
      \PP^3_{-} \ar@{-->}[r]^{\delta'_{-}} \ar[d]_{\gamma_{-}} 
      & \check{\PP}^3_{-} \ar[d]^{\check{\gamma}_{-}} \\
      \PP^8 \ar@{-->}[r]^{\calD'} & \check{\PP}^8
    } 
  \end{equation*} 
\end{prop}

\begin{proof}
  The proof is straightforward and does not depend of the duality of the Coble hypersurfaces proved
  in Theorem \ref{thm:duality}.  By \eqref{eq:F'_{01}=-F'_{02}} and \eqref{eq:F'_{b}=-F'_{-b}}, we
  see that
  \begin{equation*}
    \calD' \circ \gamma_{-} = [0:F'_{01}:-F'_{01}:F'_{10}:F'_{11}:F'_{12}:-F'_{10}:-F'_{12}:-F'_{11}].
  \end{equation*}
  On the other hand, we compute $\partial H/\partial Z_i$ and find:
  \begin{align*}
    \frac{\partial H}{\partial Z_0} &= \sum_{(i,j) \in (\FF_3)^2} \frac{\partial F}{\partial
      X_{i,j}} \frac{\partial X_{i,j}}{\partial Z_0} = 2F'_{01}, & \frac{\partial H}{\partial Z_1}
    &= F'_{10} - F'_{20} = 2F'_{10}, \\
    \frac{\partial H}{\partial Z_2} &= F'_{11} - F'_{22} = 2F'_{11}, & \frac{\partial H}{\partial Z_3}
    &= F'_{12} - F'_{21} = 2F'_{12}.
  \end{align*}
  It is then easy to check that the diagram commutes.
\end{proof}

Therefore the map $\delta'_{-}$ can be seen as the restriction of the dual map $\calD'$.

\begin{proof}[Proof of Theorem \ref{thm:dualhexahedron}]
  The restricted dual map $\delta'_{-}$ associates to a smooth point of the hexahedron the tangent
  plane through it.  But each plane is its own tangent plane so the planes contract to 6 distinct
  points in the dual $\PP^3$.  Moreover, each plane lies in the Coble sextic $\calC_6$ and is
  contracted, so the six points must be in the singular locus of Coble cubic $\calC_3$ restricted to
  $\check{\PP}^3_{-}$, which consists exactly of the 6 odd theta characteristics of our original
  curve $C$. The double cover of the unique twisted cubic through these 6 points and branched at
  these 6 points is then a curve of genus~2 having these 6 points as Weierstra\ss{} points (or odd
  theta characteristics).  So it must be isomorphic to $C$.
\end{proof}

\vspace{12pt} \subsection{The Segre cubic primal and its dual}
\label{sec:segre-cubic}

In this section we study the dual map $\calD$ restricted to $\check{\PP}^4_{+}$, which we know, by
the following theorem, is not contained in the Coble cubic $\calC_3$.

\begin{prop} \label{prop:P4 not subset of C3}
  A $(\FF_3)^4$-invariant cubic polynomial on $\check{\CC}^9$ does not vanish on
  $\check{\CC}^5_{+}$, the vector subspace whose projectivization is $\check{\PP}^4_{+}$.
\end{prop}

\begin{proof}
  The proof is easily done by looking at the equations of the $(\FF_3)^4$-invariant cubics, which,
  for some constants $\alpha_0$, \dots, $\alpha_4$, are exactly of the form \eqref{eq:C3}.  So
  restricted to $\check{\CC}^5_{+}$, as given by $\check{\gamma}_{+}$ just like in
  \eqref{eq:embeddings}, we get:
  \begin{multline*}
    \frac{\alpha_0}{3} (Y_0^3+2Y_1^3+2Y_2^3+2Y_3^3+2Y_4^3) + 2\alpha_1 (Y_0Y_1^2+2Y_2Y_3Y_4) \\
    2\alpha_2 (Y_0Y_2^2+2Y_1Y_3Y_4) + 2\alpha_3 (Y_0Y_3^2+2Y_1Y_2Y_4) + 2\alpha_4 (Y_0Y_4^2+2Y_1Y_2Y_3).
  \end{multline*}
  Clearly, the five terms are linearly independent.
\end{proof}

Thus the intersection $\check{\PP^4}_{+} \cap \calC_3$ is a cubic threefold with 10 nodes (the
maximal number possible).  It is classically known that this is only one such threefold in
$\check{\PP}^4_{+}$: the Segre cubic.  We write $\calS_3 = \check{\PP}^4_{+} \cap \calC_3$.
Following the notations of \eqref{eq:embeddings} and Section \ref{sec:det target spaces}, let
$Y_0,\dots,Y_4$ be the homogeneous coordinates of $\check{\PP}^4_{+}$ and $G$ be a homogeneous
polynomial defining the Coble cubic in $\check{\PP}^8$, so that
\begin{equation*}
  S(Y_0,Y_1,Y_2,Y_3,Y_4) = G \circ \check{\gamma}_{+} = G(Y_0,Y_1,Y_1,Y_2,Y_3,Y_4,Y_2,Y_4,Y_3)
\end{equation*}
is an equation for $\calS_3 \subset \check{\PP}^4_{+}$.  Then we have the small dual map
\begin{equation*}
  \delta_{+}:\ \check{\PP}^4_{+} \dashrightarrow \PP^4_{+}\ ,\ [Y_0:\dots:Y_4] \mapsto 
  \left[ \frac{\partial S}{\partial Y_0} : \dots : \frac{\partial S}{\partial Y_4} \right].
\end{equation*}
Similarly to Proposition \ref{prop:dualscommutewithembeddings}, we have a commutative diagram
\begin{equation} \label{eq:dualscommutewithembeddings2}
  \xymatrix{
    \check{\PP}^4_{+} \ar@{-->}[r]^{\delta_{+}} \ar[d]_{\check{\gamma}_{+}} 
    & \PP^4_{+} \ar[d]^{\gamma_{+}} \\
    \check{\PP}^8 \ar@{-->}[r]^{\calD} & \PP^8.
  } 
\end{equation}
So under the global dual map $\calD$, the Segre cubic $\calS_3$ is sent to its dual variety, the
Igusa-Segre quartic $\calI_4$, hence justifying Theorem \ref{thm:VNR} and the presence of $\calI_4$.

\vspace{12pt} \subsection{The Weddle quartic}
\label{sec:weddle-quartic}

In this section we study the dual map $\calD$ restricted to $\check{\PP}^3_{-}$ which is contained
in the Coble cubic $\calC_3$.  

Using the explicit expression for the Coble cubic explicitly given in \eqref{eq:C3}, we can see that
the image of $\check{\PP}^3_{-}$ is a $\PP^3$ of which a set of equations is
\begin{gather*}
  \alpha_0 X_{00} + \alpha_1 X_{01} + \alpha_2 X_{10} + \alpha_3 X_{12} + \alpha_4 X_{11} = 0, \\
  X_{01} - X_{02} = 0, \quad  X_{10} - X_{20} = 0, \quad X_{11} - X_{22} = 0, \quad X_{12} - X_{21} = 0.
\end{gather*}
In particular, this $\PP^3$ is contained in $\PP^4_{+}$, which confirms Proposition \ref{prop:map to
  which fixed loci}.

Moreover, this restricted dual map, which we denote $\check{\delta}_{-}$, is given by the complete
linear system (easy to check) of quadrics in $\check{\PP}^3_{-}$ passing through $J^1\cap \calC_3$,
i.e. the six odd theta characteristics of $C$.  Now the story is classically well known and we will
not prove the following result.

\begin{prop} \label{prop:w4} The rational map $\check{\delta}_{-}$ is generically 2 to 1.  The
  ramification locus is a quartic surface $\mathcal{W}_4$ with nodes at the 6 points of
  $J^1\cap\calC_3$, and the branch locus a quartic Kummer surface $\mathcal{BK}_4$.
\end{prop}

\begin{dfn}
  The ramification locus $\mathcal{W}_4$ of Proposition \ref{prop:w4} is called the \emph{\mbox{Weddle
    quartic} surface}.
\end{dfn}

This surface was first identified by T. Weddle in 1850 (\cite{Wed50} footnote p.~69).  A nice treatise
can be found in \cite{Hud05} or \cite{Cob29}.  

An important remark is that $\PP^3=\text{Im}(\check{\delta}_{-})$ lies not only in $\PP^4_{+}$ but
also in the Coble sextic $\calC_6$ because the source space $\check{\PP}^3_{-}$ lies in the Coble
cubic $\calC_3$.  From Theorem \ref{thm:VNR}, it follows that $\PP^3$ is the tangent hyperplane
$V_0$.  It is then natural to guess that the Kummer surface $\mathcal{BK}_4$ of Proposition
\ref{prop:w4} is the same as $\mathcal{K}'$ \eqref{eq:K'}, since it sits already in $V_0$.  More
generally, the 6 points in $\check{\PP}^3_{-}$ correspond to a unique curve of genus two (here it is
our original $C$), and the Kummer surface associated to the Weddle surface defined by the 6 points
should be the Kummer surface of $J(C)$.  Of course this is true (see \cite{Hud05,Bak22,Cob29}),
however, we give here a proof along the lines of the original results but putting forth previously
exhibited vector bundles of rank 3.

To identify the two Kummer surfaces, we will identify the nodes and since fixing nodes determines a
unique Kummer surface \cite{Hud05}, we get the desired result.  For this, we adapt A. Ortega
Ortega's description of the dual map on the secant variety of the embedded $J^1$ \cite[Section
7]{Ort05}.  The dual map $\calD$ is constant on secant of $J^1$, i.e. it contracts the secant.  So
the image is determined by the 2 points on $J^1$ defining the secant and there is a rational map
\begin{equation*}
  \psi:\ J^1 \times J^1 \dasharrow |3\Theta|\cong\PP^8
\end{equation*}
defining the dual map $\calD$ in the following sense.  Let $z\in|3\Theta|^*$ be a point on the
secant line passing through $a\in J^1$ and $b\in J^1$ but distinct from $a$ and $b$.  Then
\begin{equation*}
  \calD(z) = \psi(a,b).
\end{equation*}
We first define $\psi$.  Given a point $x$ on the Jacobian $J$ of the original curve $C$, we denote
by $\Theta_x$ the translate $t^*_x \Theta$ of the canonical Riemann theta divisor of $J^1$.  If we
fix two distinct points $a,b\in J^1$, then there exist $x,y \in J$ such that $\Theta_x \cap \Theta_y
= \{a,b\}$, and we set
\begin{equation*}
  \psi(a,b)= \Theta_x + \Theta_y + \Theta_{-x-y} \in |3\Theta|.
\end{equation*}

We now extend Ortega Ortega's result by showing that the divisor $\Theta_x+\Theta_y+\Theta_{-x-y}$
correspond to the vector bundle $x \oplus y \oplus (x^{-1}\otimes y^{-1})$.  Indeed, under the double
cover map $\Phi_3$ 
\eqref{eq:Phi}, the vector bundle $x
\oplus y \oplus (x^{-1}\otimes y^{-1})$ goes to the divisor
\begin{align*}
  D_{x \oplus y \oplus (x^{-1}\otimes y^{-1})} &= \{ L \in J^1 \st 
  h^0(C,(L\otimes x)\oplus(L\otimes y)\oplus(L\otimes x^{-1}\otimes y^{-1}))>0 \}, \\
  &= \Theta_x + \Theta_y + \Theta_{-x-y}.
\end{align*}
So we have proved:

\begin{lem} \label{lem:dualsecant}
  The dual map $\calD$ is defined on secants of $J^1$ by the rational map
  \begin{align*}
    \psi:\ J^1 \times J^1 &\dashrightarrow  \Sing(\SU_C(3)) \\
    (a,b)  &\longmapsto x \oplus y \oplus (x^{-1}\otimes y^{-1})
  \end{align*} 
  where $x,y \in J$ are such that $\Theta_x \cap \Theta_y = \{a,b\}$.
\end{lem}

\begin{prop} \label{prop:BK4=K'}
  The branch locus $\mathcal{BK}_4$ of the map rational map $\check{\delta}_{-}$
  is the Kummer surface $\mathcal{K}' = \{ \OO_C \oplus x \oplus x^{-1} \st x \in J \}$, i.e. the
  intersection of $V_0$ with $\calI_4$.
\end{prop}

\begin{proof}
  Let us apply Lemma \ref{lem:dualsecant}.  We know classically that the 15 secants joining any two
  points of $P$, the set of odd theta characteristics, are contracted to nodes of the Kummer surface
  $\mathcal{BK}_4$.  Let $a,b \in J^1$ be two odd theta characteristics.  We form the 2-torsion
  point $\epsilon = a-b \in J_2$, and denote by $p_{\epsilon}$ the node of $\mathcal{BK}_4$
  resulting from the contraction of the chord joining $a$ and $b$.  By definition,
  \begin{equation*}
    \Theta \cap \Theta_{\epsilon} = \{ L \in J^1 \st h^0(C,L) >0 
    \text{ and } h^0(C,L\otimes \epsilon) >0 \}.
  \end{equation*}
  We now check that $a$ and $b$ are in the intersection, therefore proving that the nodes
  $p_{\epsilon}$ correspond to vector bundles of the form $\OO_C \oplus \epsilon \oplus
  \epsilon$, which are nodes of $\mathcal{K}'$.  Indeed,
  \begin{align*}
    h^0(C,\OO_C(a))&> 0 & h^0(C,\OO_C(a+a-b)) = h^0(C,\OO_C(b)) &> 0, \\
    h^0(C,\OO_C(b)) &> 0 & h^0(C,\OO_C(b+a-b)) = h^0(C,\OO_C(a)) &> 0.
  \end{align*}
  Thus we have proved that the Kummer surfaces $\mathcal{BK}_4$ and $\mathcal{K}'$ share 15 nodes.
  Their well-known configuration guarantees that there is a conic on which 6 of the 15 nodes lie.
  The double cover of this curve, branched over the 6 points, is a curve of genus two, whose
  Jacobian is the unique (up to isomorphism) abelian surface covering the Kummer surfaces.  So the
  Kummer surfaces $\mathcal{BK}_4$ and $\mathcal{K}'$ are isomorphic and the curve of genus two is
  isomorphic to our original curve $C$, since $\mathcal{K}'$ is the Kummer surface associated to
  $C$.  Moreover, $\mathcal{BK}_4$ and $\mathcal{K}'$ have the same tropes because of the way we
  identify the 15 nodes.  Since each pair of tropes intersects at exactly 2 points, actually nodes,
  there has to be some pair of tropes that meet at only one node of the previously identified 15
  nodes, in which case the second point of intersection is the last and sixteenth node.  Thus, the
  two Kummer surfaces are actually equal.
\end{proof}

This proposition can be restated as follows: in the projective space $\PP^3$, six points in general
position define a curve $C$ of genus two.  The complete linear system of quadric surfaces through
those 6 points defines a rational double cover of another $\PP^3$ whose branch locus is the Kummer
quartic surface of $C$.

As an immediate corollary of the proof, we obtain:

\begin{cor}
  The twisted cubic curve through the 6 theta characteristics of $\check{\PP}^3_{-}$ is contracted
  under the dual map $\calD$ to the point of tangency of $V_0$ with the Igusa-Segre quartic
  $\calI_4$.
\end{cor}





\vspace{18pt}
\bibliography{vbcgg2}

\providecommand{\bysame}{\leavevmode\hbox to3em{\hrulefill}\thinspace}
\providecommand{\MR}{\relax\ifhmode\unskip\space\fi MR }
\providecommand{\MRhref}[2]{%
  \href{http://www.ams.org/mathscinet-getitem?mr=#1}{#2}
}
\providecommand{\href}[2]{#2}
\begin{thebibliography}{BNR89}

\bibitem[Ati57]{Ati57b}
M.~F. Atiyah, \emph{Vector bundles over an elliptic curve}, Proc. London Math.
  Soc. (3) \textbf{7} (1957), 414--452.

\bibitem[Bak22]{Bak22}
H.~Baker, \emph{Principles of geometry}, Cambridge Univ. Press, Cambridge,
  1922.

\bibitem[Bal90]{Bal90}
V.~Balaji, \emph{Intermediate jacobians of some moduli spaces of vector bundles
  on curves}, Amer. J. Math. \textbf{112} (1990), 611--630.

\bibitem[Bar95]{Bar95}
W.~Barth, \emph{Quadratic equations for level-3 abelian surfaces}, in ``Abelian
  varieties'' (Egloffstein, 1993), de Gruyter, Berlin, 1995, pp.~1--18.

\bibitem[Bea03]{Bea03}
A.~Beauville, \emph{The {Coble} hypersurfaces}, C. R. Acad. Sci. Paris, Ser. I
  \textbf{336} (2003), no.~3, 189--194.

\bibitem[BNR89]{BNR89}
A.~Beauville, M.~S. Narasimhan, and S.~Ramanan, \emph{Spectral curves and the
  generalised theta divisor}, J. reine angew. Math. \textbf{398} (1989),
  169--179.

\bibitem[Cob17]{Cob17}
A.~Coble, \emph{Point sets and allied {C}remona groups {III}}, Trans. Amer.
  Math. Soc. \textbf{18} (1917), 331--372.

\bibitem[Cob61]{Cob29}
\bysame, \emph{Algebraic geometry and theta functions}, American Mathematical
  Society Colloquium Publication, vol.~X, Amer. Math. Soc., Providence, RI,
  1961, Revised reprint of the 1929 original.

\bibitem[DN89]{DN89}
J.~M. Drezet and M.~S. Narasimhan, \emph{Groupe de {P}icard des vari\'et\'es de
  modules de fibr\'es semi-stables sur les courbes alg\'ebriques}, Invent.
  Math. \textbf{97} (1989), 53--94.

\bibitem[EL92]{EL92}
L.~Ein and R.~Lazarsfeld, \emph{Stability and restrictions of {P}icard bundles,
  with an application to the normal bundles of elliptic curves}, in ``Complex
  Projective Geometry'', Cambridge University Press, 1992, pp.~149--156.

\bibitem[FH79]{FH79}
W.~Fulton and J.~Hansen, \emph{A connectedness theorem for projective varities,
  with applications to intersections and singularities of mappings}, Ann. of
  Math. \textbf{110} (1979), no.~1, 159--166.

\bibitem[Gee82]{vdG82}
G.~van~der Geer, \emph{On the geometry of a {S}iegel modular threefold}, Math.
  Ann. \textbf{260} (1982), no.~3, 317--350.

\bibitem[GKZ94]{GKZ94}
I.~M. Gelfand, M.~M. Kapranov, and A.~V. Zelevinski, \emph{Discriminants,
  resultants, and multidimensional detreminants}, Mathematics: {T}heories \&
  {A}pplications, Birkhauser Boston, Inc., Boston, MA, 1994.

\bibitem[GLP83]{GLP83}
L.~Gruson, R.~Lazarsfeld, and C.~Peskine, \emph{On a theorem of {C}astelnuovo,
  and the equations defining space curves}, Invent. Math. \textbf{72} (1983),
  no.~3, 491--506.

\bibitem[Gro68]{SGA2}
A.~Grothendieck, \emph{Cohomologie locale des faisceaux coh\'erents et
  th\'eor\`emes de {L}efschetz locaux et globaux {$(SGA\ 2)$}}, North-Holland
  Publishing Co., Amsterdam; Masson \& Cie, Editeur, Paris, 1968.

\bibitem[HR04]{HR04}
J.-M. Hwang and S.~Ramanan, \emph{Hecke curves and {H}itchin discriminant},
  Ann. Sci. Ecole Norm. Sup. (4) \textbf{37} (2004), no.~5, 801--817.

\bibitem[Hud90]{Hud05}
R.~W. H.~T. Hudson, \emph{Kummer's quartic surface}, Cambridge University
  Press, Cambridge, 1990, Revised reprint of the 1905 original.

\bibitem[Koi76]{Koi76}
S.~Koizumi, \emph{Theta relations and projective normality of abelian
  varieties}, Amer. J. Math. \textbf{98} (1976), 865--889.

\bibitem[KP95]{KP95}
A.~Kouvidakis and T.~Pantev, \emph{The automorphism group of the moduli space
  of semi stable vector bundles}, Math. Ann. \textbf{302} (1995), no.~2,
  225--268.

\bibitem[Las96]{Las96}
Y.~Laszlo, \emph{Local structure of the moduli space of vector bundles over
  curves}, Comment. Math. Helv. \textbf{71} (1996), no.~3, 373--401.

\bibitem[LB92]{LB92}
H.~Lange and C.~Birkenhake, \emph{Complex abelian varieties}, Springer-Verlag,
  1992.

\bibitem[MN68]{MN68}
D.~Mumford and P.~E. Newstead, \emph{Periods of a moduli space of bundles on
  curves}, Amer. J. Math. \textbf{90} (1968), 1201--1208.

\bibitem[New68]{New68}
P.~E. Newstead, \emph{Stable bundles of rank 2 and odd degree over a curve of
  genus 2}, Topology \textbf{7} (1968), 205--215.

\bibitem[Ngu05]{Ngu05}
Q.~M. Nguy{$\tilde{\hat{e}}$}n, \emph{Dualities and classical geometry of the
  moduli space of rank 3 vector bundle on a curve of genus 2}, {Ph.D}
  dissertation, University of Michigan, 2005.

\bibitem[NR69]{NR69}
M.~S. Narasimhan and S.~Ramanan, \emph{Moduli of vector bundles on a compact
  {R}iemann surface}, Ann. of Math. \textbf{89} (1969), 14--51.

\bibitem[NR75]{NR75}
\bysame, \emph{Deformations of the moduli space of vector bundles on a curve},
  Ann. of Math. \textbf{101} (1975), 391--417.

\bibitem[NR87]{NR87}
\bysame, \emph{$2\theta$-linear system on abelian varieties}, Vector bundles on
  algebraic varieties (Bombay 1984), Tata Inst. Fund. Res. Stud. Math.,
  vol.~11, Tata Inst. Fund. Res., Bombay, 1987, pp.~415--427.

\bibitem[NR03]{NgRa03}
Q.~M. Nguy{$\tilde{\hat{e}}$}n and S.~Rams, \emph{On the geometry of the
  {C}oble-{D}olgachev sextic}, Le Matematiche (Catania) \textbf{58} (2003),
  no.~2, 257--275.

\bibitem[OO05]{Ort05}
A.~Ortega~Ortega, \emph{On the moduli space of rank 3 vector bundles on a genus
  2 curve and the {Coble} cubic}, J. Algebraic Geom. \textbf{14} (2005), no.~2,
  327--356.

\bibitem[Pau02]{Pau02}
C.~Pauly, \emph{Self-duality of {Coble's} quartic hypersurface and
  application}, Michigan Math. J. \textbf{50} (2002), 551--574.

\bibitem[Ses67]{Ses67}
C.~S. Seshadri, \emph{Space of unitary vector bundles on a compact {R}iemann
  surface}, Ann. of Math. \textbf{85} (1967), 303--336.

\bibitem[Tu93]{Tu93}
L.~W. Tu, \emph{Semistable bundles over an elliptic curve}, Adv. Math.
  \textbf{98} (1993), 1--26.

\bibitem[Tyu74]{Tyu74}
A.~N. Tyurin, \emph{The geometry of moduli of vector bundles}, Russ. Math.
  Surveys \textbf{29} (1974), no.~6, 57--88.

\bibitem[Wed50]{Wed50}
T.~Weddle, \emph{On the theorems in space analogous to those of {Pascal} and
  {Brianchon} in a plane.---part~{II}}, Camb. and Dub. Math. Jour. \textbf{5}
  (1850), 58--69.

\end{thebibliography}

\end{document}